\documentclass[12pt]{amsart}
\usepackage{amsmath,amssymb,amsthm}
\usepackage[numbers,sort&compress]{natbib}
\usepackage{enumitem}
\usepackage{geometry}
\usepackage{tikz-cd}
\usepackage[colorlinks=true,linkcolor=blue,citecolor=blue,urlcolor=blue]{hyperref}

\newcommand{\Z}{\mathbb Z}
\newcommand{\B}{\mathfrak B}
\newcommand{\X}{X(G_k)}
\newcommand{\Lk}{\operatorname{Lk}}
\newcommand{\Span}{\operatorname{span}}
\newcommand{\Hom}{\operatorname{Hom}}
\newcommand{\rank}{\operatorname{rank}}

\title[Homological description of $\mathcal Z_{k+1}(G_k)$]
{Homological Description of Equivariant Geometric Bordism}
\author{Bo Chen, Hao Li and Zhi L\"u}

\address{School of Mathematics and Statistics, Huazhong University of Science and Technology, Wuhan, China}
\email{bobchen@hust.edu.cn}
\address{College of Mathematics and Sciences, Shanghai Normal University\newline
Shanghai 200234, China}
\email{14110840001@fudan.edu.cn}
\address{School of Mathematical Sciences, Fudan University, Shanghai, 200433, P. R. China}
\email{zlu@fudan.edu.cn}

\thanks{Partially supported by NSFC grants 11971112 and 12501087.}

\theoremstyle{plain}
\newtheorem{theorem}{Theorem}[section]
\newtheorem{corollary}[theorem]{Corollary}
\newtheorem{lemma}[theorem]{Lemma}
\newtheorem{proposition}[theorem]{Proposition}
\theoremstyle{definition}
\newtheorem{definition}[theorem]{Definition}
\newtheorem{notation}[theorem]{Notation}

\newtheorem{example}[theorem]{Example}
\newtheorem{remark}[theorem]{Remark}

\allowdisplaybreaks
\begin{document}

\begin{abstract}

Using the evenness criterion in \cite{LCLS}  and duality principle in \cite{CLT}, we construct  a finite chain complex
$\mathfrak B_\bullet((\mathbb Z_2)^k)$ built from the universal complex
$X((\mathbb Z_2)^k)$ and its links. We then show that the equivariant unoriented bordism group $\mathcal Z_{k+1}((\mathbb Z_2)^k)$ of all $(k+1)$-dimensional smooth  closed connected manifolds with effective $\Z_2^k$-actions fixing isolated points,  is naturally isomorphic to 
$H_{k-2}(\mathfrak B_\bullet((\mathbb Z_2)^k);\mathbb Z_2)$.  The vertical
spectral sequence associated with the natural double-complex structure
on $\mathfrak B_\bullet((\mathbb Z_2)^k)$ collapses at $E^3$, yielding
an explicit 
formula for $\dim_{\mathbb Z_2}\mathcal Z_{k+1}((\mathbb Z_2)^k)$ for
every $k$.
\end{abstract}

\maketitle

\section{Introduction}\label{sec:introduction}

Throughout this paper, $G_k=(\mathbb Z_2)^k$ is the elementary abelian $2$-group of rank $k$. Conner and Floyd \cite{CF} initiated the study of
the equivariant geometric bordism ring $\mathcal Z_*(G_k)=\sum_n\mathcal{Z}_n(G_k)$ represented by smooth  closed connected manifolds with effective $\Z_2^k$-actions fixing isolated points.  They established the algebraic framework by introducing
the fixed-point
homomorphism
\[
 \phi_*\colon\mathcal Z_*(G_k)\longrightarrow\mathcal R_*(G_k),
 \qquad
 \phi_*([M])=\sum_{p\in M^{G_k}}[T_pM].
\]
Here $\mathcal R_*(G_k)$ is the polynomial algebra over $\mathbb Z_2$
generated by the isomorphism classes of irreducible real
$G_k$-representations.  Every such representation is one-dimensional
and is determined by a character
$\rho\in V_k:=\Hom(G_k,\mathbb Z_2)$, which we call a \emph{weight}.
Thus an $n$-dimensional representation is recorded by the formal
product of its weights, called its
\emph{weight monomial}.
Conner and Floyd determined the ring structure of $\mathcal Z_*(G_k)$
for $k=1,2$ \cite[Theorems~25.1, 31.1 and 31.2]{CF}.  For $k\geq3$,
the full ring structure remains open.  Partial results include studies
of indecomposable elements \cite{M,MS} and the recent proof that
$\mathcal Z_*(G_k)$ is not finitely generated \cite{GLN}.

 A classical theorem of Stong
\cite{St} shows that $\phi_*$ is injective.  We may therefore
identify $\mathcal Z_n(G_k)$ with $\operatorname{im}\phi_n$ and ask
which sums of faithful degree-$n$ weight monomials lie in this image.

The first general answer to this realization problem came from
localization.  tom Dieck \cite{d} proved that equivariant
characteristic numbers characterize geometric fixed-point data by an
integrality condition.  Kosniowski and Stong \cite{KS} subsequently
expressed the relevant equivariant characteristic numbers in terms of fixed-point
data by a localization formula.  The resulting tom
Dieck--Kosniowski--Stong criterion (see \cite[Theorem 1.1]{LCLS} for details) characterizes realizability, for
all $n$ and $k$, by the integrality of localization sums.  While
complete as a theoretical characterization, its direct application
can be difficult: the localization formulas involve denominators and
all symmetric polynomials, and do not immediately provide either a
finite test stated on the given weight monomials or a structural
description of $\operatorname{im}\phi_n$.

Nevertheless, using Davis--Januszkiewicz theory \cite{DJ}, Chen, L\"u,
and Tan \cite{CLT} obtained a complete answer for $n=k$, the smallest
degree in which a faithful weight monomial for $G_k$ exists.  In this
dimension, the weights of a faithful monomial form a basis of $V_k$.
Using the differential introduced by L\"u and Tan \cite{LT} on the
dual of the Conner--Floyd representation algebra, they viewed the
dual basis as a maximal simplex of the universal complex $X(G_k)$ and
identified $\mathcal Z_k(G_k)$ with its top homology.
Here $X(G_k)$ is the simplicial complex whose vertices are the
nonzero elements of $G_k$ and whose simplices are their $\Z_2$-linearly
independent subsets.  This gives both a global description of all
realizable data and an effective computation of the group, thereby
resolving both issues in the case $n=k$. 

The companion paper \cite[Theorem~1.2]{LCLS} replaces the localization
criterion by a finite denominator-free \emph{evenness condition},
stated directly in terms of weight multiplicities after grouping the
weight monomials with respect to a nonzero $\rho\in V_k$ and
restriction to $\ker\rho$.
It has yielded complete calculations of the groups in two low-degree cases
$(n,k)=(4,3)$ and $(5,3)$, including both their dimensions and
geometric generators \cite{LCLS,GL}.

Two difficulties nevertheless remain.  First, the evenness condition is
an elementwise membership test, not a global algebraic description of
$\operatorname{im}\phi_n$.  Second, as $n$ and $k$ increase, the
number of weights, equivalence classes, and parity equations grows
rapidly.  Although the test can be implemented on a computer, a direct
calculation quickly becomes unwieldy and gives little structural
insight.

Combining the evenness theorem with the dual viewpoint developed above, we obtain a complete computation of $\mathcal Z_{k+1}(G_k)$ for every $k$ in this paper.

Let $\mathcal A$ be a finite set of faithful degree-$(k+1)$ weight
monomials for $G_k$ with no zero weight factor.
For each nonzero weight $\rho$, the evenness condition
\cite[Definition~2.4]{LCLS} first partitions the monomials containing
$\rho$ into classes $\mathcal A_{\rho^d,j}$, where $d$ is the
multiplicity of $\rho$.  In degree $k+1$, the multiplicity bound
\cite[Remark~1]{LCLS} gives $d\leq2$.  The evenness theorem
\cite[Theorem~1.2]{LCLS} says that $\mathcal A$ is realizable as
fixed-point data if and only if, for every $0\ne\rho\in V_k$ and every
such class,
\[
\begin{array}{ll}
 d=1:&|\mathcal A_{\rho^1,j}|\equiv0\pmod2,\\[3pt]
 d=2:&|\mathcal A_{\rho^2,j}|
       \equiv m(\mathcal A_{\rho^2,j},\beta)
       \equiv0\pmod2
       \quad\text{for every }0\ne\beta\in V_k,
\end{array}
\]
where $m(\mathcal C,\beta)$ denotes the total multiplicity, among the
monomials in $\mathcal C$, of the weight $\beta$.
Equivalently, these conditions characterize when
$\sum_{\tau\in\mathcal A}\tau$ belongs to
$\operatorname{im}\phi_{k+1}$.

The first line is a parity test on the size of each class.  Once a
weight is repeated, the second line requires, in addition, a parity
test for every weight occurring across that class.  These tests must
be carried out for every nonzero $\rho$, so a direct application of
the evenness theorem produces a large family of parity equations and
gives no efficient way to determine their joint solution space.

The multiplicity-two conditions admit a decisive dual reformulation.
For $\tau\in\mathcal A_{\rho^2,j}$, write
$\tau=\rho^2\rho_2\cdots\rho_k$.  Then
$\{\rho,\rho_2,\ldots,\rho_k\}$ is a basis of
$V_k$.  Let $\alpha_\tau(\rho)\in G_k$ be the basis
vector dual to $\rho$, called the \emph{dual pivot}.
Proposition~\ref{prop:dual-pivot} shows that all parity equations in
the second line for $\mathcal A_{\rho^2,j}$ are equivalent to the
single vector equation
\[
 \sum_{\tau\in\mathcal A_{\rho^2,j}}
 \alpha_\tau(\rho)=0.
\]
Thus the weight-by-weight parity tests for a multiplicity-two class
reduce to the vanishing of one vector sum.

Our first step is to linearize and assemble all these classwise tests;
its key ingredient is the dual-pivot criterion above.  The linear map
$\Omega_0$ encodes the first line: each multiplicity-one class
contributes its class monomial with coefficient equal to the
cardinality of the class (Proposition~\ref{prop:Omega0}).  The map
$\Omega_1$ assembles the vector
equations for the second line: each multiplicity-two class determines
a simplex $\sigma_0(\rho,j)$, and tensoring its dual-pivot sum with
this simplex keeps the contributions of distinct classes separate
(Proposition~\ref{prop:Omega1}).  Together they encode the entire
evenness condition by the single equation recorded in
Corollary~\ref{cor:global-obstruction}:
\[
 \Omega\left(\sum_{\tau\in\mathcal A}\tau\right)=0,
 \qquad
 \Omega=(\Omega_0,\Omega_1).
\]

Our second step, motivated by the Chen--L\"u--Tan duality in the case
$n=k$, turns this kernel into homology.  For a faithful
degree-$(k+1)$ weight monomial, the weights involved in the unique
dependence dualize to an equivalence chain $c(\sigma_1)$, while the
remaining independent weights
dualize to a complementary simplex $\sigma_0$ in the common link of the
simplices occurring in $c(\sigma_1)$.  Their tensor product is independent of the chosen
basis.  In this way we obtain a linear isomorphism $D$ from the space
spanned by such weight monomials onto the top chain group of a finite complex
$\mathfrak B_\bullet(G_k)$; see Proposition~\ref{prop:d-iso} and
\eqref{eq:B-k-2}.  Applying the corresponding dualization maps $D$ to
the source and target of $\Omega$ identifies $\Omega$ with the
boundary $\partial_{k-2}$, as made explicit in
Definition~\ref{def:boundary-op} and
Proposition~\ref{prop:boundary-obstruction}.  Consequently,
Theorem~\ref{thm:cycle-condition} gives, for every finite set
$\mathcal A$ of such weight monomials,
\[
 \sum_{\tau\in\mathcal A}\tau\in\operatorname{im}\phi_{k+1}
 \quad\Longleftrightarrow\quad
 \partial_{k-2}
 \left(\sum_{\tau\in\mathcal A}D(\tau)\right)=0.
\]

This cycle criterion turns the elementwise evenness theorem into a
homological description of the entire group.
\begin{theorem}\label{thm:homology-isomorphism}
For every $k\geq1$, dualization induces an isomorphism
\[
 \mathcal Z_{k+1}(G_k)
 \cong
 H_{k-2}\bigl(\B_\bullet(G_k)\bigr).
\]
\end{theorem}

The complex $\B_\bullet(G_k)$ has a natural double-complex structure whose nonexceptional vertical
columns are augmented simplicial chain complexes of links in $X(G_k)$.
The link topology forces the vertical spectral sequence to collapse at $E^3$.

Let $f_{p,k}$ be the number of $p$-simplices of
$X(G_k)$ and let $A_{p,k}$ be the number of
top-dimensional spheres in the link of a $p$-simplex, with
$A_{-1,k}$ referring to $X(G_k)$ itself.  Counting the surviving
terms in the spectral sequence gives the following closed formula.
\begin{theorem}\label{thm:intro-dimension}
For $k\geq2$,
\[
\dim_{\mathbb Z_2}\mathcal Z_{k+1}(G_k)
=A_{0,k}f_{0,k}
 +\sum_{p=1}^{k-2}\frac{A_{p,k}f_{p,k}}{p+2}
 +kA_{-1,k}
 -\left(k-\frac{1}{k+1}\right)f_{k-1,k}.
\]
\end{theorem}

The quantities in this formula are given explicitly by
\eqref{eq:a-k}, \eqref{eq:a-pk}, and \eqref{eq:f-pk}, so the formula
is directly computable.  The value for $k=1$, which lies outside the
range of Theorem~\ref{thm:intro-dimension}, follows from a direct
computation of the spectral sequence in this case.  Together, the
first values are
\[
\begin{array}{|c|c|c|c|c|c|c|c|}
\hline
k&1&2&3&4&5&6&7\\ \hline
\dim_{\mathbb Z_2}\mathcal Z_{k+1}(G_k)
&0&0&32&3{,}177&719{,}164&476{,}303{,}715
&1{,}025{,}099{,}895{,}814\\
\hline
\end{array}
\]
 The two low-rank values
$\mathcal Z_2(G_1)=0$ and $\mathcal Z_3(G_2)=0$ agree with the
classical calculations of Conner and Floyd
\cite[Theorems~25.1, 31.1, and~31.2]{CF}.
In particular, for $k=3$ the formula gives
$\dim_{\mathbb Z_2}\mathcal Z_4(G_3)=32$, in agreement with the
calculation in the companion paper \cite{LCLS}.

Section~\ref{sec:background} states the degree-$(k+1)$ evenness criterion.
Section~\ref{sec:dual-construction} constructs the dual data and the dualization isomorphism.
Section~\ref{sec:dual-data} derives the cycle criterion.
Section~\ref{sec:chain-complex} constructs the double complex and computes the homology of its total complex via the vertical spectral sequence.
The appendix verifies the horizontal-differential identities.

\section{Fixed point data and the evenness criterion}
\label{sec:background}
We first recall the general notation for weight monomials and mod
$\rho$ equivalence from the companion paper \cite[Section~2]{LCLS},
and then state the degree-$(k+1)$ evenness theorem used below.  All
vector spaces and chain groups in this paper are over $\mathbb Z_2$.

\begin{notation}[Weight monomials]\label{not:weight-monomials}
We use $G_k=(\mathbb Z_2)^k$ both as an elementary abelian $2$-group
and as its underlying $k$-dimensional vector space over $\mathbb Z_2$.
Put
\[
 V_k=\Hom(G_k,\mathbb Z_2),
 \qquad V_k^\times=V_k\setminus\{0\}.
\]
Every irreducible real $G_k$-representation is one-dimensional and is
determined by a character $\rho\in V_k$; we refer to such a $\rho$ as
a \emph{weight}.  We encode an $n$-dimensional real
$G_k$-representation by the formal product of its weights
\[
 \tau=\rho_1\cdots\rho_n,
 \qquad \rho_i\in V_k,
\]
which is its \emph{weight monomial}.  Its degree is $\deg\tau=n$.
The trivial character $\mathbf0\in V_k$ is the \emph{zero weight} and
records a trivial one-dimensional summand.
The Conner--Floyd representation algebra is the polynomial algebra
\[
 \mathcal R_*(G_k)
 \cong
 \mathbb Z_2[x_\rho \mid \rho\in V_k],
\]
where $x_\rho$ represents the one-dimensional representation with
weight $\rho$.  Under this identification, the notation
$\rho_1\cdots\rho_n$ abbreviates
$x_{\rho_1}\cdots x_{\rho_n}$.  In particular, a zero weight factor
$\mathbf0$ denotes the variable $x_{\mathbf0}$, not the zero element
of the polynomial algebra; hence a weight monomial containing
$\mathbf0$ is nonzero.
As in
\cite[Section~2.1]{LCLS}, the monomial
determines the representation up to isomorphism, so we work
exclusively with weight monomials.  We denote the multiplicity of
$\rho\in V_k$ in $\tau$ by
$m(\tau,\rho)$, following \cite[Definition~2.2]{LCLS}.  If
$\mathcal C$ is a finite set of weight monomials for $G_k$, set, for
$\rho\in V_k$,
\[
 m(\mathcal C,\rho)=\sum_{\tau\in\mathcal C}m(\tau,\rho).
\]
Only the parity of these integers will be used.
\end{notation}

\begin{notation}[Deletion and restriction]
\label{not:deletion-restriction}
Fix $\rho\in V_k^\times$.  If
\[
 \tau=\rho^d\rho_1\cdots\rho_{n-d},
 \qquad \rho_i\ne\rho,
\]
then
\begin{equation}\label{eq:delete-all-restrict}
 (\tau/\rho^d)|_{\ker\rho}
 :=(\rho_1|_{\ker\rho})\cdots
   (\rho_{n-d}|_{\ker\rho}).
\end{equation}
Thus one first deletes all $d$ occurrences of $\rho$ and then
restricts every remaining weight to $\ker\rho$.  If $\tau$ has no
zero weight factor, then every restricted weight in
\eqref{eq:delete-all-restrict} is nonzero, since $\rho$ is the only
nonzero weight vanishing on $\ker\rho$.

We will also use the related notation
\begin{equation}\label{eq:delete-one-restrict}
 (\tau/\rho)|_{\ker\rho},
\end{equation}
in which only one occurrence of $\rho$ is deleted before restriction.
If $m(\tau,\rho)=2$, the undeleted copy of $\rho$ restricts to the
zero character
\[
 \mathbf0\in\Hom(\ker\rho,\mathbb Z_2),
\]
which corresponds to the trivial one-dimensional representation of
$\ker\rho$.  We retain this zero-weight factor in
$(\tau/\rho)|_{\ker\rho}$; this convention will be used when
dualizing \eqref{eq:delete-one-restrict} in
Subsection~\ref{subsec:deleted-restriction}.
\end{notation}

\begin{definition}[Mod $\rho$ equivalence classes]
\label{def:mod-rho-classes}
Fix $\rho\in V_k^\times$.  Following
\cite[Definition~2.3]{LCLS}, two weight monomials $\tau$ and $\tau'$
for $G_k$ are
\emph{mod $\rho$ equivalent}, written $\tau\sim_\rho\tau'$, if, for
some $d\geq1$,
\[
 m(\tau,\rho)=m(\tau',\rho)=d,
 \qquad
 (\tau/\rho^d)|_{\ker\rho}
 =(\tau'/\rho^d)|_{\ker\rho}.
\]
Let $\mathcal A$ be a finite set of degree-$n$ weight monomials for
$G_k$, and set
\[
 \mathcal A_\rho
 :=\{\tau\in\mathcal A\mid\rho\mid\tau\}.
\]
For each $d\geq1$, let $J_{\rho,d}$ index the mod $\rho$
equivalence classes in $\mathcal A_\rho$ whose elements have
multiplicity $d$ in $\rho$, with the sets $J_{\rho,d}$ chosen
pairwise disjoint as $d$ varies.  Denoting these classes by
$\mathcal A_{\rho^d,j}$, we have
\[
 \mathcal A_\rho
 =\bigsqcup_{d\geq1}\ \bigsqcup_{j\in J_{\rho,d}}
   \mathcal A_{\rho^d,j}.
\]
We call $\mathcal A_{\rho^d,j}$ a
\emph{multiplicity-$d$ mod $\rho$ equivalence class}.
\end{definition}

We now specialize to the degree relevant to this paper.  Recall that
a weight monomial for $G_k$ is \emph{faithful} if its weight factors
span $V_k$.  In a faithful degree-$(k+1)$ weight monomial, every
nonzero weight has multiplicity at most two: if $d\geq3$, the
number of distinct weights is at most
$1+(k+1-d)\leq k-1$; so its weights cannot span $V_k$.  This is the specialization of the
multiplicity bound in \cite[Remark~1]{LCLS}.

Combining the realization theorem \cite[Theorem~1.2]{LCLS} with the
evenness condition of \cite[Definition~2.4]{LCLS} and the preceding
multiplicity bound gives exactly the following statement.

\begin{theorem}[Evenness theorem in degree $k+1$]
\label{thm:evenness}
Let $\mathcal A$ be a finite set of faithful degree-$(k+1)$ weight
monomials for $G_k$ with no zero weight factor.  Then
\[
 \sum_{\tau\in\mathcal A}\tau\in\operatorname{im}\phi_{k+1}
\]
if and only if, for every $\rho\in V_k^\times$,
\begin{equation}\label{eq:evenness-one}
 |\mathcal A_{\rho^1,j}|\equiv0\pmod2
 \qquad \forall j\in J_{\rho,1},
\end{equation}
and 
\begin{equation}\label{eq:evenness-two}
   |\mathcal A_{\rho^2,j}| \equiv m(\mathcal{A}_{\rho^2,j},\beta) \equiv 0 \pmod{2}
  \qquad \forall \beta\in V_k^\times, j\in J_{\rho,2}.
\end{equation}
\end{theorem}

We call \eqref{eq:evenness-one} and \eqref{eq:evenness-two} the
\emph{multiplicity-one condition} and the \emph{multiplicity-two
condition}, respectively.  Together they form the \emph{evenness
condition} in degree $k+1$.

\begin{notation}[Class monomials]\label{not:class-monomials}
For a nonempty class
$\mathcal A_{\rho^d,j}$, the definition of mod $\rho$ equivalence
implies that the monomial
$ (\tau/\rho^d)|_{\ker\rho}$ is independent of the choice of
$\tau\in\mathcal A_{\rho^d,j}$.  We may therefore write
\begin{equation}\label{eq:class-monomial}
 g_{\rho,j}:=(\tau/\rho^d)|_{\ker\rho}
 \qquad(\tau\in\mathcal A_{\rho^d,j}).
\end{equation}
We call $g_{\rho,j}$ the \emph{class monomial} of
$\mathcal A_{\rho^d,j}$.
Thus $g_{\rho,j}$ is a degree-$(k+1-d)$ weight monomial for
$\ker\rho$.  We regard it as a $\ker\rho$-labelled weight monomial;
the label records that its factors belong to
$\Hom(\ker\rho,\mathbb Z_2)$.
\end{notation}

\begin{lemma}[Separation of class monomials]
\label{lem:class-monomials}
Let $\mathcal A$ be a finite set of degree-$(k+1)$ weight monomials
for $G_k$ with no zero weight factor.  For the nonempty
classes $\mathcal A_{\rho^d,j}$ of
Definition~\ref{def:mod-rho-classes}, the assignment
\[
 \mathcal A_{\rho^d,j}
 \longmapsto g_{\rho,j}
\]
is injective as $\rho$, $d$, and $j$ vary.
\end{lemma}

\begin{proof}
For fixed $\rho$, the degree of $g_{\rho,j}$ determines $d$, since
$\deg g_{\rho,j}=k+1-d$.  Once $d$ is fixed, equality of two such
monomials is precisely the condition that their representatives are
mod $\rho$ equivalent, and hence that the two classes coincide.  If
$\rho\ne\rho'$, then $\ker\rho\ne\ker\rho'$: over $\mathbb Z_2$, a
nonzero linear function is uniquely determined by its kernel.
Hence $g_{\rho,j}$ and $g_{\rho',j'}$ have different labels and are
therefore distinct as labelled weight monomials.
\end{proof}


\section{Dual data and the dualization isomorphism}
\label{sec:dual-construction}

Throughout this section, let $G\subseteq G_k$ be a subgroup of rank
$r$, regarded  also as an $r$-dimensional $\mathbb Z_2$-vector space.  Its
\emph{universal complex} $X(G)$ has vertex set $G\setminus\{0\}$ and
simplices the linearly independent subsets.  We use the augmented
simplicial chains $\widetilde C_*(X(G))$, with
$\widetilde C_{-1}(X(G))=\mathbb Z_2\langle\emptyset\rangle$ and
$\widetilde C_q(X(G))=0$ for $q<-1$.

The terminology of Notation~\ref{not:weight-monomials} applies with
$G_k$ replaced by $G$.  A weight monomial for $G$ is faithful if its
weight factors span $\Hom(G,\mathbb Z_2)$.

\subsection{Equivalence chains in the universal complex}
\label{subsec:equiv-rel}

\begin{definition}\label{def:equiv-sim}
For simplices $\sigma,\sigma'\in X(G)$, write $\sigma\sim\sigma'$ if
$\sigma=\sigma'$, or if they share a vertex $\alpha$ and
\[
 \sigma'\setminus\{\alpha\}
 = (\sigma\setminus\{\alpha\})+\alpha,
 \qquad
 S+\alpha:=\{\beta+\alpha\mid\beta\in S\}.
\]
In either case, we call $\sigma$ and $\sigma'$ \emph{equivalent}.
\end{definition}

\begin{proposition}\label{prop:equiv-facts}
The relation $\sim$ has the following properties.
\begin{enumerate}
 \item Equivalent simplices have the same cardinality; two distinct
 equivalent simplices have exactly one common vertex.
 \item It is an equivalence relation on the simplices of $X(G)$.
 \item The empty simplex and every $0$-simplex form singleton classes.
 \item The equivalence class of a $p$-simplex has $p+2$ elements
 whenever $0<p<r$.
 \item Equivalent simplices span the same subspace of $G$ and have
 the same link in $X(G)$.
\end{enumerate}
\end{proposition}

\begin{proof}
The empty simplex and every $0$-simplex are equivalent only to
themselves.  For
$\sigma=\{\alpha_1,\dots,\alpha_{p+1}\}$ with $0<p<r$, the
simplices related to $\sigma$ are $\sigma$ and
\[
 \{\alpha_i\}\cup
 \{\alpha_j+\alpha_i\mid j\ne i\},
 \qquad 1\le i\le p+1.
\]
These $p+2$ simplices are distinct.  Any two shifted simplices are
related through their common vertex $\alpha_i+\alpha_j$, so they form
one equivalence class.  This proves (1)--(4), including transitivity.
Each shift is a change
of basis in $\Span_{\mathbb Z_2}\sigma$, and a link in a universal
complex depends only on this span, proving (5).
\end{proof}

\begin{notation}[Equivalence chains]
For a $p$-simplex $\sigma\in X(G)$, let $\overline\sigma$ denote its
equivalence class and set
\begin{equation}\label{eq:equiv-sum}
 c(\sigma):=\sum_{\sigma'\in\overline\sigma}\sigma'
 \in\widetilde C_p(X(G)).
\end{equation}
Thus $c(\emptyset)=\emptyset$, $c(\{\alpha\})=\{\alpha\}$, and, for
$p>0$,
\[
 c(\sigma)=\sigma+
 \sum_{\alpha\in\sigma}
 \bigl((\sigma\setminus\{\alpha\})+\alpha\bigr)\cup\{\alpha\}.
\]
\end{notation}

\begin{lemma}[Shift equivalence]\label{lem:shift-equivalence}
Let $\sigma\in X(G)$ have positive dimension.  If
$\alpha,\beta\in\sigma$ are distinct, then
\[
 (\sigma\setminus\{\alpha\})+\alpha
 \sim
 (\sigma\setminus\{\beta\})+\beta.
\]
\end{lemma}

\begin{proof}
The two simplices share $\alpha+\beta$, and translation by this
vertex identifies their remaining vertices.
\end{proof}

\subsection{Dual data in rank \texorpdfstring{$r$}{r}}
\label{subsec:dual-data}

Let $\tau$ be a faithful degree-$(r+1)$ weight monomial for
$G$.  Choose a basis
$\{\rho_1,\dots,\rho_r\}$ of
$\Hom(G,\mathbb Z_2)$ among its factors and write
\[
 \tau=\rho_1\cdots\rho_r\gamma,
 \qquad
 \gamma=\sum_{i\in\Lambda}\rho_i.
\]
If $\gamma=\mathbf0$, then $\Lambda=\emptyset$.  The
evaluation isomorphism $G\cong\Hom(G,\mathbb Z_2)^*$ identifies the
dual basis of $\{\rho_1,\dots,\rho_r\}$ with a basis
$\{\alpha_1,\dots,\alpha_r\}$ of $G$.  Set
\[
 \sigma_1=\{\alpha_i\mid i\in\Lambda\},
 \qquad
 \sigma_0=\{\alpha_i\mid i\notin\Lambda\}.
\]

\begin{definition}\label{def:dual-data}
With the notation above, the \emph{dual data} of $\tau$ is
\begin{equation}\label{eq:d-def}
 D(\tau):=c(\sigma_1)\otimes\sigma_0
 \in \Span_{\mathbb Z_2}\{c(\sigma_1)\}\otimes
 \widetilde C_{r-|\Lambda|-1}
 \bigl(\Lk_{X(G)}(\sigma_1)\bigr).
\end{equation}
\end{definition}

\begin{remark}\label{rem:dual-data-link}
Since $\sigma_1\sqcup\sigma_0$ is a basis of $G$, we have
$\sigma_0\in\Lk_{X(G)}(\sigma_1)$.  Proposition
~\ref{prop:equiv-facts}(5) shows that this link is common to all
simplices occurring in $c(\sigma_1)$.  The space
$\Span_{\mathbb Z_2}\{c(\sigma_1)\}$ is the one-dimensional subspace
of $\widetilde C_{|\Lambda|-1}(X(G))$ generated by this equivalence
chain.
\end{remark}

\begin{lemma}\label{lem:dual-data-well-defined}
Let $\tau$ be a faithful degree-$(r+1)$ weight monomial for
$G$.  Then $D(\tau)$ is independent of the basis  chosen
among the factors of $\tau$.
\end{lemma}

\begin{proof}
If the remaining factor is $\gamma=\mathbf0$, the other $r$ factors
form the unique basis among the factors of $\tau$.
Suppose $\gamma\ne\mathbf0$.  Every other basis among the factors is
obtained by replacing some $\rho_{i_0}$, with $i_0\in\Lambda$, by
$\gamma$; all such replacements are bases because the coefficient of
$\rho_{i_0}$ in $\gamma=\sum_{i\in\Lambda}\rho_i$ is one.  Let
$\{\alpha_1',\dots,\alpha_r'\}$ be the basis of $G$ dual to the new
weight basis, with $\alpha_{i_0}'$ paired with $\gamma$.  Directly
from the evaluation pairing,
\[
 \alpha_{i_0}'=\alpha_{i_0},
 \qquad
 \alpha_i'=
 \begin{cases}
  \alpha_i+\alpha_{i_0},&i\in\Lambda\setminus\{i_0\},\\
  \alpha_i,&i\notin\Lambda.
 \end{cases}
\]
The remaining weight is $\rho_{i_0}=\gamma+\sum_{ i\in \Lambda \setminus\{i_0\}} \rho_i$, whose dependence on the new
basis gives the partition
\[
 \sigma_1'
 =\{\alpha_{i_0}\}\cup
   ((\sigma_1\setminus\{\alpha_{i_0}\})+\alpha_{i_0}),
 \qquad
 \sigma_0'=\sigma_0.
\]
Thus $\sigma_1'\sim\sigma_1$, so
$c(\sigma_1')=c(\sigma_1)$. Since $\sigma_0'=\sigma_0$, the resulting
dual data is unchanged.
\end{proof}

\begin{definition}[Weight-monomial spaces]\label{def:F-spaces}
Let $\mathcal F(G)$ be the $\mathbb Z_2$-vector space whose basis is
the set of all faithful degree-$(r+1)$ weight monomials for
$G$.  Let $\mathcal F_0(G)\subseteq\mathcal F(G)$ be the subspace
whose basis consists of those monomials having no zero weight factor.
\end{definition}

\begin{definition}\label{def:B-space}
With $G$ and $r$ as fixed above, for integers $p,q$ set
\begin{equation}\label{eq:B-pq}
 B_{p,q}(G):=
 \bigoplus_{\overline{\sigma^p}}
 \Span_{\mathbb Z_2}\{c(\sigma^p)\}\otimes
 \widetilde C_q\bigl(\Lk_{X(G)}(\sigma^p)\bigr),
\end{equation}
where the sum runs over equivalence classes of $p$-simplices, and put
\[
 \widehat\B_{r-2}(G):=\bigoplus_{p+q=r-2}B_{p,q}(G).
\]
\end{definition}

By Remark~\ref{rem:dual-data-link}, each summand in \eqref{eq:B-pq}
is well-defined independently of the representative of
$\overline{\sigma^p}$.  The sum is direct because equivalence chains
from distinct classes have disjoint supports.

\begin{proposition}\label{prop:d-iso}
The assignment $\tau\mapsto D(\tau)$ extends linearly to an
isomorphism
\[
 D\colon\mathcal F(G)\overset{\cong}{\longrightarrow}
 \widehat\B_{r-2}(G).
\]
\end{proposition}

\begin{proof}
For a basis element
$c(\sigma_1)\otimes\sigma_0\in B_{p,r-2-p}(G)$, write
\[
 \sigma_1=\{\alpha_1,\dots,\alpha_{p+1}\},
 \qquad
 \sigma_0=\{\alpha_{p+2},\dots,\alpha_r\},
\]
and let $\rho_1,\dots,\rho_r$ be the dual basis.  With
$\rho_0=\sum_{i=1}^{p+1}\rho_i$, define
\begin{equation}\label{eq:d-inverse}
 D^{-1}(c(\sigma_1)\otimes\sigma_0)
 :=\rho_0\rho_1\cdots\rho_r.
\end{equation}
For $p=-1$, we use the conventions $\sigma_1=\emptyset$ and
$\rho_0=\sum_{i=1}^{0}\rho_i=\mathbf0$.

If $p>0$, replacing $\sigma_1$ by the representative
$\{\alpha_1,\alpha_1+\alpha_2,\dots,
\alpha_1+\alpha_{p+1}\}$, after relabelling if necessary, replaces
the relevant part of the dual
basis by $\rho_0,\rho_2,\dots,\rho_{p+1}$.  The corresponding
dependent weight is
$\rho_0+\sum_{i=2}^{p+1}\rho_i=\rho_1$, and hence
\[
 \rho_1\rho_0\rho_2\cdots\rho_r
 =\rho_0\rho_1\cdots\rho_r.
\]
The cases $p=-1,0$ have singleton equivalence classes.  Since every
representative of $\overline{\sigma_1}$ is obtained by one of the
displayed replacements, $D^{-1}$ is well-defined.

Now let $\tau=\rho_0\rho_1\cdots\rho_r$, where
$\{\rho_1,\dots,\rho_r\}$ is a basis and
$\rho_0=\sum_{i=1}^{p+1}\rho_i$.  Its dual data is
$c(\sigma_1)\otimes\sigma_0$, so \eqref{eq:d-inverse} gives
$D^{-1}D(\tau)=\tau$.  Conversely, starting with
$c(\sigma_1)\otimes\sigma_0$, formula \eqref{eq:d-inverse} produces
exactly this basis presentation, and Definition~\ref{def:dual-data}
returns the original tensor.  Hence $DD^{-1}$ is also the identity.
\end{proof}

\begin{remark}\label{rem:b-components}
Proposition~\ref{prop:d-iso} identifies the components of
$\widehat\B_{r-2}(G)$ with the following types of degree-$(r+1)$ weight monomials.  The
component
\[
 B_{-1,r-1}(G)
 =\mathbb Z_2\langle\emptyset\rangle\otimes
  \widetilde C_{r-1}(X(G))
\]
corresponds to monomials having a zero weight factor.  For
$0\le p\le r-1$, a basis element of $B_{p,r-2-p}(G)$ is
$c(\sigma_1)\otimes\sigma_0$, where
$\sigma_1\sqcup\sigma_0$ is a basis of $G$; the unique circuit of
the corresponding weight monomial has length $p+2$.   Moreover,
$B_{p,q}(G)=0$ for $p<-1$ or $p\ge r$. Thus
$\widehat\B_{r-2}(G)$ assembles the dual data of all faithful
degree-$(r+1)$ weight monomials for $G$. 
The hat distinguishes this full dual-data space from the groups below
associated  with monomials having no zero weight factor.
\end{remark}

\subsection{Specialization to \texorpdfstring{$G_k$}{Gk} and its
hyperplanes}
\label{subsec:deleted-restriction}

Specializing Proposition~\ref{prop:d-iso} to $G=G_k$ and to
monomials without a zero weight factor gives
\begin{equation}\label{eq:B-k-2}
 D\colon\mathcal F_0(G_k)\overset{\cong}{\longrightarrow}
 \B_{k-2}(G_k),
 \qquad
 \B_{k-2}(G_k)
 :=\bigoplus_{p=0}^{k-1}B_{p,k-p-2}(G_k).
\end{equation}
Here and below, the spaces $B_{p,q}(G_k)$ are those defined in
\eqref{eq:B-pq}.
For $\rho\in V_k^\times$, the subgroup $\ker\rho$ has rank $k-1$,
and the same construction gives
\begin{equation}\label{eq:d-kernel}
 D\colon\mathcal F(\ker\rho)
 \overset{\cong}{\longrightarrow}\widehat\B_{k-3}(\ker\rho).
\end{equation}

Using the delete-one--restrict convention of
Notation~\ref{not:deletion-restriction}, for $\tau$ and $\rho$ as in
the proposition below we have
\[
 (\tau/\rho)|_{\ker\rho}\in\mathcal F(\ker\rho).
\]
The following formulas relate its dual data to $D(\tau)$ and will be
used to define the boundary operator in
Section~\ref{sec:dual-data}.

\begin{proposition}\label{prop:d-deleted-restriction}
Let $\tau$ be a faithful degree-$(k+1)$ weight monomial for $G_k$
with no zero weight factor.  Choose a basis
$\mathcal B=\{\rho_1,\dots,\rho_k\}$ of $V_k$ among the factors of
$\tau$, with dual basis
$\{\alpha_1,\dots,\alpha_k\}$, and write
$D(\tau)=c(\sigma_1)\otimes\sigma_0$.  For a weight $\rho$ occurring
as a factor of $\tau$, the following hold.
\begin{enumerate}[label=\textup{(\roman*)}]
 \item If $\rho\in\mathcal B$ and $\alpha$ is its dual vertex, then
 \[
 D((\tau/\rho)|_{\ker\rho})=
 \begin{cases}
  c(\sigma_1\setminus\{\alpha\})\otimes\sigma_0,
    &\alpha\in\sigma_1,\\[4pt]
  c(\sigma_1)\otimes(\sigma_0\setminus\{\alpha\}),
    &\alpha\in\sigma_0.
 \end{cases}
 \]
 \item If $\rho\notin\mathcal B$, write
 $\rho=\sum_{i\in\Lambda}\rho_i$, so that
 $\sigma_1=\{\alpha_i\mid i\in\Lambda\}$.  For any
 $i_0\in\Lambda$,
 \[
 D((\tau/\rho)|_{\ker\rho})
 =c((\sigma_1\setminus\{\alpha_{i_0}\})+\alpha_{i_0})
  \otimes\sigma_0,
 \]
 independently of the choice of $i_0$.
\end{enumerate}
\end{proposition}

\begin{proof}
In case \textup{(i)}, suppose after relabelling that $\rho=\rho_1$.
Restriction induces an isomorphism
\[
 V_k/\langle\rho_1\rangle
 \cong\Hom(\ker\rho_1,\mathbb Z_2),
\]
so $\rho_2|_{\ker\rho_1},\dots,
\rho_k|_{\ker\rho_1}$ form a basis whose dual basis is
$\alpha_2,\dots,\alpha_k$.  Thus deleting $\rho_1$ deletes its dual
vertex $\alpha_1$ from the dual-basis partition.  According as
$\alpha_1$ belongs to $\sigma_1$ or $\sigma_0$, this gives the first
or the second displayed formula.

In case \textup{(ii)}, the factor outside $\mathcal B$ is the
dependent weight
$\rho=\sum_{i\in\Lambda}\rho_i$.  Hence its restrictions satisfy
$\sum_{i\in\Lambda}\rho_i|_{\ker\rho}=0$.  For any
$i_0\in\Lambda$, the weights
$\{\rho_i|_{\ker\rho}\mid i\ne i_0\}$ form a basis of
$\Hom(\ker\rho,\mathbb Z_2)$.  Its dual basis is
\[
 \{\alpha_i+\alpha_{i_0}\mid
     i\in\Lambda\setminus\{i_0\}\}
 \sqcup
 \{\alpha_i\mid i\notin\Lambda\},
\]
so the two components of the dual data are
$(\sigma_1\setminus\{\alpha_{i_0}\})+\alpha_{i_0}$ and
$\sigma_0$, respectively.  This gives the stated formula, whose
independence of $i_0$ follows from
Lemma~\ref{lem:shift-equivalence}.
\end{proof}

\begin{example}\label{ex:repeated-weight}
    Let $\tau = \rho_1^2 \rho_2 \cdots \rho_k$, so $\sigma_1 = \{\alpha_1\}$ and $\sigma_0 = \{\alpha_2,\dots,\alpha_k\}$. Then
    \[
    \begin{aligned}
    D(\tau) &= c(\{\alpha_1\})\otimes \{\alpha_2, \dots, \alpha_k\} = \{\alpha_1\}\otimes \{\alpha_2, \dots, \alpha_k\}, \\
    D((\tau/\rho_1)|_{\ker\rho_1}) &= \emptyset \otimes \{\alpha_2, \dots, \alpha_k\}, \\
    D((\tau/\rho_i)|_{\ker\rho_i}) &= \{\alpha_1\}\otimes \{\alpha_2, \dots, \hat{\alpha}_i, \dots, \alpha_k\}, \quad i = 2,\dots,k.
    \end{aligned}
    \]
    Here $\rho_1$ participates in the dependence ($\alpha_1 \in \sigma_1$), so its removal shrinks $\sigma_1$ to $\emptyset$. The other weights $\rho_i$ ($i \ge 2$) lie in $\sigma_0$, so their removal deletes $\alpha_i$ from $\sigma_0$.
\end{example}

\begin{example}\label{ex:boundary}
    Let $\tau = \rho_0\rho_1\rho_2\cdots\rho_k$, where $\{\rho_1,\dots,\rho_k\}$ is a basis and $\rho_0 = \rho_1 + \rho_2$. Then $\sigma_1 = \{\alpha_1,\alpha_2\}$ and $\sigma_0 = \{\alpha_3,\dots,\alpha_k\}$. By Proposition~\ref{prop:d-deleted-restriction},
    \[
    \begin{aligned}
    D(\tau) &= c(\{\alpha_1,\alpha_2\})\otimes \{\alpha_3, \dots, \alpha_k\}, \\
    D((\tau/\rho_0)|_{\ker\rho_0}) &= c(\{\alpha_{12}\}) \otimes \{\alpha_3, \dots, \alpha_k\}, \\
    D((\tau/\rho_1)|_{\ker\rho_1}) &= c(\{\alpha_2\}) \otimes \{\alpha_3, \dots, \alpha_k\}, \\
    D((\tau/\rho_2)|_{\ker\rho_2}) &= c(\{\alpha_1\}) \otimes \{\alpha_3, \dots, \alpha_k\}, \\
    D((\tau/\rho_i)|_{\ker\rho_i}) &= c(\{\alpha_1,\alpha_2\})\otimes \{\alpha_3, \dots, \hat{\alpha}_i, \dots, \alpha_k\}, \quad i = 3,\dots,k.
    \end{aligned}
    \]
\end{example}

\section{Evenness criterion as a cycle condition}\label{sec:dual-data}
Let $\mathcal A$ be a finite set of faithful degree-$(k+1)$ weight
monomials for $G_k$ with no zero weight factor.
We now convert the degree-$(k+1)$ evenness criterion on $\mathcal A$ into the
vanishing of one boundary operator.  The construction is organized
in two stages.  First we encode multiplicity-one and multiplicity-two
mod $\rho$ equivalence classes by two global obstruction maps
$\Omega_0$ and $\Omega_1$.  We then transport their direct sum
through the dualization isomorphisms of the preceding section.

\subsection{\texorpdfstring
{The zeroth obstruction: multiplicity-one mod $\rho$ classes}
{The zeroth obstruction: multiplicity-one classes}}
By Definition~\ref{def:F-spaces}, $\mathcal F_0(G)$ has as basis the
faithful degree-$(r+1)$ weight monomials for $G$ with no zero weight
factor, where $r=\operatorname{rank}G$.  For a nonzero weight
$\rho\in V_k^\times$, we
have $\ker\rho\cong G_{k-1}$.  Define
\begin{equation}\label{eq:Omega0}
\begin{aligned}
 \Omega_0\colon\mathcal F_0(G_k)
 &\longrightarrow
 \bigoplus_{\rho\in V_k^\times}\mathcal F_0(\ker\rho),\\
 \Omega_0(\tau)
 &=
 \sum_{\substack{\rho\mid\tau\\m(\tau,\rho)=1}}
 (\tau/\rho)|_{\ker\rho}.
\end{aligned}
\end{equation}

\begin{proposition}\label{prop:Omega0}
For $\mathcal A$ as above,
\[
 \Omega_0\left(\sum_{\tau\in\mathcal A}\tau\right)=0
\]
if and only if, for every $\rho\in V_k^\times$, each multiplicity-one mod
$\rho$ equivalence class in $\mathcal A$ satisfies the
multiplicity-one condition \eqref{eq:evenness-one}.
\end{proposition}

\begin{proof}
By \eqref{eq:class-monomial}, a multiplicity-one mod $\rho$
equivalence class $\mathcal A_{\rho^1,j}$ determines the unique basis
element $g_{\rho,j}=(\tau/\rho)|_{\ker\rho}$ of $\mathcal F_0(\ker\rho)$, independent of the choice of $\tau\in \mathcal A_{\rho^1,j}$.  Hence
\[
 \Omega_0\left(\sum_{\tau\in\mathcal A}\tau\right)
 =\sum_{\rho,j}|\mathcal A_{\rho^1,j}|\,g_{\rho,j}.
\]
By Lemma~\ref{lem:class-monomials}, these class monomials are distinct
basis elements in the direct sum indexed by $\rho$.  The assertion
follows.
\end{proof}

\subsection{The first obstruction: the dual-pivot moment}

Fix $\rho\in V_k^\times$.  Let $\mathcal A_{\rho^2,j}$ be a
multiplicity-two mod $\rho$ equivalence class, and take
\[
 \tau=\rho^2\rho_2\cdots\rho_k
 \in\mathcal A_{\rho^2,j}.
\]
Up to ordering,
$\{\rho,\rho_2,\ldots,\rho_k\}$ is the unique basis of $V_k$ that
can be chosen among its factors.  Let
\[
 \{\alpha_\tau(\rho),\alpha_\tau(\rho_2),\dots,
   \alpha_\tau(\rho_k)\}
\]
be its dual basis in $G_k$.  We call
$\alpha_\tau(\rho)$ the \emph{dual pivot} of the pair
$(\tau,\rho)$ and set
\begin{equation}\label{eq:kernel-simplex}
 \sigma_0(\tau,\rho)
 :=\{\alpha_\tau(\rho_2),\dots,\alpha_\tau(\rho_k)\}
 \in X(\ker\rho).
\end{equation}
Indeed, the vertices of this simplex form the basis of $\ker\rho$
dual to
$\rho_2|_{\ker\rho},\dots,\rho_k|_{\ker\rho}$.
By mod $\rho$ equivalence, this restricted basis, and hence its dual
basis, is independent of $\tau\in\mathcal A_{\rho^2,j}$.  We may
therefore write
\begin{equation}\label{eq:class-simplex}
 \sigma_0(\rho,j):=\sigma_0(\tau,\rho)
 \in X(\ker\rho).
\end{equation}

\begin{proposition}[Dual-pivot criterion for repeated weights]
\label{prop:dual-pivot}
Let $\mathcal A_{\rho^2,j}$ be a multiplicity-two mod $\rho$
equivalence class.  Then it satisfies the multiplicity-two condition
\eqref{eq:evenness-two} if and only if
\[
 \sum_{\tau\in\mathcal A_{\rho^2,j}}
 \alpha_\tau(\rho)=0
 \quad\text{in }G_k.
\]
\end{proposition}

\begin{proof}
Choose
\[
 \tau_0=\rho^2\rho_2\cdots\rho_k
\]
in the class and let
$(\alpha_{\tau_0}(\rho),\alpha_2,\ldots,\alpha_k)$ be the dual
basis of $(\rho,\rho_2,\ldots,\rho_k)$.  Every other member $\tau=\rho^2\rho_2'\cdots\rho_k'$ has
\[
 \rho_i'=\rho_i+\epsilon_i(\tau)\rho,
 \qquad \epsilon_i(\tau)\in\mathbb F_2, ~i=2,\cdots, k,
\]
and its dual pivot is
\[
 \alpha_\tau(\rho)
 =
 \alpha_{\tau_0}(\rho)
 +\sum_{i=2}^k\epsilon_i(\tau)\alpha_i.
\]
Therefore
\[
 \sum_{\tau\in\mathcal A_{\rho^2,j}}\alpha_\tau(\rho)
 =
 |\mathcal A_{\rho^2,j}|\alpha_{\tau_0}(\rho)
 +\sum_{i=2}^k
 \left(\sum_\tau\epsilon_i(\tau)\right)\alpha_i.
\]
This sum vanishes if and only if the class has even cardinality and
$\sum_\tau\epsilon_i(\tau)=0$ for every $i\geq2$.

For each $i\geq2$, exactly one of $\rho_i$ and $\rho_i+\rho$ occurs
in $\tau$, and
\[
 m(\tau,\rho_i+\rho)=\epsilon_i(\tau),
 \qquad
 m(\tau,\rho_i)=1+\epsilon_i(\tau)
\]
in $\mathbb F_2$.  Thus the displayed conditions are equivalent to 
\[
 |\mathcal A_{\rho^2,j}|
 \equiv m(\mathcal A_{\rho^2,j},\rho_i)
 \equiv m(\mathcal A_{\rho^2,j},\rho_i+\rho)
 \equiv0\pmod2
 \qquad(2\le i\le k).
\]
The weight $\rho$ itself occurs twice in every monomial, and all
weights outside these pairs occur zero times.  Hence these
congruences are exactly the multiplicity-two condition
\eqref{eq:evenness-two} for $\mathcal A_{\rho^2,j}$.
\end{proof}

Define the global first-moment obstruction
\begin{equation}\label{eq:Omega1}
\begin{aligned}
 \Omega_1\colon\mathcal F_0(G_k)
 &\longrightarrow
 G_k\otimes\widetilde C_{k-2}(\X),\\
 \Omega_1(\tau)
 &=
 \sum_{\substack{\rho\mid\tau\\m(\tau,\rho)=2}}
 \alpha_\tau(\rho)\otimes\sigma_0(\tau,\rho).
\end{aligned}
\end{equation}

\begin{proposition}\label{prop:Omega1}
For $\mathcal A$ as above,
\[
 \Omega_1\left(\sum_{\tau\in\mathcal A}\tau\right)=0
\]
if and only if, for every $\rho\in V_k^\times$, each multiplicity-two mod
$\rho$ equivalence class in $\mathcal A$ satisfies the
multiplicity-two condition \eqref{eq:evenness-two}.
\end{proposition}

\begin{proof}
By \eqref{eq:class-simplex}, each nonempty multiplicity-two mod
$\rho$ equivalence class $\mathcal A_{\rho^2,j}$ determines the
simplex
$\sigma_0(\rho,j)=\sigma_0(\tau,\rho)$ for any
$\tau\in\mathcal A_{\rho^2,j}$.  The assignment
\[
 (\rho,j)\longmapsto\sigma_0(\rho,j)
\]
is injective.  Indeed, if
$\sigma_0(\rho,j)=\sigma_0(\rho',j')$, then taking spans gives
$\ker\rho=\ker\rho'$, and hence $\rho=\rho'$ over $\mathbb Z_2$.
The common simplex determines a unique dual weight monomial for this
kernel $\ker \rho$, so $g_{\rho,j}=g_{\rho,j'}$;
Lemma~\ref{lem:class-monomials} now gives $j=j'$.

Consequently, the simplices $\sigma_0(\rho,j)$ indexing these classes
are distinct, and
\[
 \Omega_1\left(\sum_{\tau\in\mathcal A}\tau\right)
 =
 \sum_{\mathcal A_{\rho^2,j}\ne\varnothing}
 \left(
   \sum_{\tau\in\mathcal A_{\rho^2,j}}
   \alpha_\tau(\rho)
 \right)
 \otimes\sigma_0(\rho,j).
\]
These simplices are distinct basis elements of
$\widetilde C_{k-2}(\X)$ and are therefore linearly independent.
Thus the left-hand side vanishes if and only if every inner vector
sum vanishes, which by
Proposition~\ref{prop:dual-pivot} is equivalent to the asserted
multiplicity-two condition.
\end{proof}

Put
\[
 \Omega=(\Omega_0,\Omega_1).
\]

\begin{corollary}[Global obstruction form of evenness]
\label{cor:global-obstruction}
For $\mathcal A$ as above,
\[
 \mathcal A\text{ satisfies the evenness condition}
 \quad\Longleftrightarrow\quad
 \Omega\left(\sum_{\tau\in\mathcal A}\tau\right)=0.
\]
\end{corollary}

\begin{proof}
By Theorem~\ref{thm:evenness}, evenness in degree $k+1$ is the
conjunction of the multiplicity-one condition
\eqref{eq:evenness-one} and the multiplicity-two condition
\eqref{eq:evenness-two}.  These are respectively equivalent to
$\Omega_0(\sum_{\tau\in\mathcal A}\tau)=0$ and
$\Omega_1(\sum_{\tau\in\mathcal A}\tau)=0$ by
Propositions~\ref{prop:Omega0} and~\ref{prop:Omega1}.
\end{proof}

\subsection{Dualization and the boundary operator}

We now transport the obstruction map $\Omega=(\Omega_0,\Omega_1)$
to the dual-data spaces.  Its source is dualized by
\eqref{eq:B-k-2}:
\[
 D\colon\mathcal F_0(G_k)\overset{\cong}{\longrightarrow}
 \B_{k-2}(G_k)
 :=\bigoplus_{p=0}^{k-1}B_{p,k-p-2}(G_k).
\]
Denote the target of $\Omega$ by
\begin{equation}\label{eq:Omega-target}
 \mathcal T_{k-3}(G_k):=
 \bigl(G_k\otimes\widetilde C_{k-2}(\X)\bigr)
 \oplus
 \bigoplus_{\rho\in V_k^\times}\mathcal F_0(\ker\rho).
\end{equation}
The first summand, which receives $\Omega_1$, is already in
dual-data form.  To dualize the summand receiving $\Omega_0$, we
restrict \eqref{eq:d-kernel} to obtain
\begin{equation}\label{eq:d-kernel-zero}
 D\colon\mathcal F_0(\ker\rho)
 \overset{\cong}{\longrightarrow}\B_{k-3}(\ker\rho),
 \qquad
 \B_{k-3}(\ker\rho)
 :=\bigoplus_{p=0}^{k-2}B_{p,k-p-3}(\ker\rho).
\end{equation}

A basis element
$c(\sigma_1)\otimes\sigma_0\in B_{p,k-p-3}(G_k)$ satisfies
$|\sigma_1|+|\sigma_0|=k-1$.  Hence
$\sigma_1\sqcup\sigma_0$ is a basis of a unique hyperplane
$H\subset G_k$, and there is a unique $\rho\in V_k^\times$ with
$H=\ker\rho$.  This gives the canonical decomposition
\begin{equation}\label{eq:hyperplane-decomposition}
 \bigoplus_{p=0}^{k-2}B_{p,k-p-3}(G_k)
 \cong
 \bigoplus_{\rho\in V_k^\times}
 \B_{k-3}(\ker\rho).
\end{equation}
Combining \eqref{eq:d-kernel-zero} and
\eqref{eq:hyperplane-decomposition}, we therefore dualize the target
\eqref{eq:Omega-target} as
\begin{equation}\label{eq:B-k-3}
 \B_{k-3}(G_k)
 :=
 \bigl(G_k\otimes\widetilde C_{k-2}(\X)\bigr)
 \oplus
 \bigoplus_{p=0}^{k-2}B_{p,k-p-3}(G_k),
\end{equation}
where the spaces $B_{p,q}(G_k)$ are those of \eqref{eq:B-pq}.
Under \eqref{eq:hyperplane-decomposition}, the identity on the first
summand together with $\bigoplus_\rho D$ identifies
$\mathcal T_{k-3}(G_k)$ with \eqref{eq:B-k-3}.

\begin{definition}\label{def:boundary-op}
The \emph{boundary operator} $\partial_{k-2}$ is the transport of
$\Omega$ under the preceding dualizations; equivalently, it is the
unique map making the following diagram commute:
\[
\begin{tikzcd}[column sep=large,row sep=large]
 \mathcal F_0(G_k)
 \arrow[r,"{(\Omega_0,\Omega_1)}"]
 \arrow[d,"D"',"\cong"]
 &
 \mathcal T_{k-3}(G_k)
 \arrow[d,"{\mathrm{id}\,\oplus\,\bigoplus_\rho D}","\cong"']
 \\
 \B_{k-2}(G_k)
 \arrow[r,"\partial_{k-2}"']
 &
 \B_{k-3}(G_k).
\end{tikzcd}
\]
\end{definition}

\begin{proposition}\label{prop:boundary-obstruction}
For
$c(\sigma_1)\otimes\sigma_0\in B_{p,k-p-2}(G_k)$, the boundary is
given as follows.  Here $d^X$ denotes the boundary map in the augmented
simplicial chain complex $\widetilde C_\bullet(\X)$.  If $p\geq1$,
choose $\alpha\in\sigma_1$.  Then
\[
\begin{aligned}
 \partial_{k-2}(c(\sigma_1)\otimes\sigma_0)
 ={}&
 \left(
 c((\sigma_1\setminus\{\alpha\})+\alpha)
 +\sum_{\beta\in\sigma_1}c(\sigma_1\setminus\{\beta\})
 \right)\otimes\sigma_0\\
 &+c(\sigma_1)\otimes d^X(\sigma_0).
\end{aligned}
\]
The first equivalence chain is independent of the choice of
$\alpha$ by Lemma~\ref{lem:shift-equivalence}.  If $p=0$ and
$\sigma_1=\{\alpha\}$, then
\[
 \partial_{k-2}(c(\{\alpha\})\otimes\sigma_0)
 =
 \alpha\otimes\sigma_0
 +c(\{\alpha\})\otimes d^X(\sigma_0).
\]
In particular, for every finite set $\mathcal A$ of faithful
degree-$(k+1)$ weight monomials for $G_k$ with no zero weight factor,
\[
 \partial_{k-2}
 \left(\sum_{\tau\in\mathcal A}D(\tau)\right)=0
 \quad\Longleftrightarrow\quad
 \Omega\left(\sum_{\tau\in\mathcal A}\tau\right)=0.
\]
\end{proposition}

\begin{proof}
If the circuit simplex $\sigma_1$ has dimension $p\geq1$, the
weight monomial has no repeated weight.  Proposition
~\ref{prop:d-deleted-restriction} shows that dualizing all multiplicity-one
deletion--restriction terms in $\Omega_0$ gives
\[
 \left(
 c((\sigma_1\setminus\{\alpha\})+\alpha)
 +\sum_{\beta\in\sigma_1}
 c(\sigma_1\setminus\{\beta\})
 \right)\otimes\sigma_0
 +c(\sigma_1)\otimes d^X(\sigma_0).
\]
If $p=0$, the weight monomial has one repeated weight.  Its
contribution to $\Omega_1$ is
$\alpha\otimes\sigma_0$, while its multiplicity-one
deletion--restriction terms
give $c(\{\alpha\})\otimes d^X(\sigma_0)$.  These are exactly the two
displayed formulas.  The final equivalence follows from the
commutative diagram in Definition~\ref{def:boundary-op}.
\end{proof}

\subsection{The cycle criterion}

\begin{theorem}\label{thm:cycle-condition}
Let $\mathcal A$ be a finite set of faithful degree-$(k+1)$
weight monomials for $G_k$  with no zero weight factor.  Then
\[
 \sum_{\tau\in\mathcal A}\tau \in\operatorname{im}\phi_{k+1}
 \quad\Longleftrightarrow\quad
 \partial_{k-2}
 \left(\sum_{\tau\in\mathcal A}D(\tau)\right)=0.
\]
\end{theorem}

\begin{proof}
By Theorem~\ref{thm:evenness} (the evenness theorem),
$\sum_{\tau\in\mathcal A}\tau\in\operatorname{im}\phi_{k+1}$  if and
only if $\mathcal A$ satisfies evenness condition.  By
Corollary~\ref{cor:global-obstruction}, this is equivalent to
$\Omega(\sum_{\tau\in\mathcal A}\tau)=0$.  Proposition
~\ref{prop:boundary-obstruction} identifies the latter condition
with the displayed cycle equation.
\end{proof}

\begin{corollary}\label{cor:top-exact}
The sequence
\[
 0\longrightarrow\mathcal Z_{k+1}(G_k)
 \xrightarrow{D\circ\phi_{k+1}}
 \B_{k-2}(G_k)
 \xrightarrow{\partial_{k-2}}
 \B_{k-3}(G_k)
\]
is exact.
\end{corollary}

\section{The dual-data complex and its homology}\label{sec:chain-complex}

The boundary operator $\partial_{k-2}$ constructed in Section~\ref{sec:dual-data} extends naturally to a chain complex. We formalize this as the total complex of a double complex $E_{p,q}$, and then compute its homology via a spectral sequence.

\subsection{\texorpdfstring
{The double complex $E_{p,q}$}
{The double complex E(p,q)}}\label{subsec:double-complex}

For $-1 \leq p \leq k-1$ and $q \geq -1$, define
\[
E_{p,q} = 
\begin{cases}
\bigoplus_{\overline{\sigma^p}} \Span\{c(\sigma^p)\} \otimes \widetilde C_q(\Lk_{\X} \sigma^p), & p \geq 0, \\[4pt]
G_k \otimes \widetilde C_q(\X), & p = -1, -1\leq q\leq k-2, \\[4pt]
0, &\text{otherwise}.
\end{cases}
\]
Set $E_{p,q}=0$ outside this range.  For $p\geq0$, we have
$E_{p,q}=B_{p,q}(G_k)$; the exceptional column $E_{-1,q}$ is modeled
on the tensor form occurring in $\operatorname{im}\Omega_1$.

We extend the two components of $\partial_{k-2}$ to all bidegrees as
follows.  We write $d^X$ for the simplicial boundary.  The
\emph{vertical differential} $(d_v)_{p,q} \colon E_{p,q} \to E_{p,q-1}$
is induced by $d^X$:
\[
(d_v)_{p,q}(c(\sigma^p) \otimes \sigma^q)
=c(\sigma^p) \otimes d^X(\sigma^q).
\]

The \emph{horizontal differential} $(d_h)_{p,q} \colon E_{p,q} \to E_{p-1,q}$ is defined by
\[
(d_h)_{p,q}(c(\sigma^p) \otimes \sigma^q) = 
\begin{cases}
\alpha \otimes \sigma^q, & p = 0,\ \sigma^p = \{\alpha\}, \\[4pt]
c(\sigma') \otimes \sigma^q + \sum_{\alpha \in \sigma^p} c(\sigma^p \setminus \{\alpha\}) \otimes \sigma^q, & p \geq 1,
\end{cases}
\]
where $\sigma' = (\sigma^p \setminus \{\alpha\}) + \alpha$ for any chosen $\alpha \in \sigma^p$ (the result is independent of the choice; see Lemma~\ref{lem:shift-equivalence}).

\begin{proposition}\label{prop:double-complex}
The maps $d_h$ and $d_v$ define a double complex
$(E_{p,q},d_h,d_v)$ over $\mathbb Z_2$.
\end{proposition}

\begin{proof}
The identity $d_v^2=0$ follows from the simplicial boundary.  The
horizontal differential is well-defined by
Lemma~\ref{lem:dh-well-defined}, and $d_h^2=0$ by
Proposition~\ref{prop:dh-square-zero}.  Finally, $d_h$ changes only
the first tensor factor, whereas $d_v$ applies the simplicial
boundary to the second; hence $d_hd_v=d_vd_h$ on every generator,
including the map into the exceptional column.
\end{proof}

\subsection{The total complex and homology isomorphism}

Let
\[
 \B_i(G_k):=\bigoplus_{p+q=i}E_{p,q},
 \qquad \partial_i:=d_h+d_v.
\]
This is the total complex of $E_{p,q}$.  Since $E_{p,q}=0$ for
$p+q\geq k-1$, its top chain group is
$\B_{k-2}(G_k)$ and
\[
 H_{k-2}(\B_\bullet(G_k))=\ker\partial_{k-2}.
\]
In this degree, the chain groups and boundary agree with those
constructed in Section~\ref{sec:dual-data}.

\begin{proof}[Proof of Theorem~\ref{thm:homology-isomorphism}]
The preceding identification and Corollary~\ref{cor:top-exact} give
the assertion.
\end{proof}

We henceforth assume $k\ge2$ and compute the homology of the total
complex by a spectral sequence.

\subsection{\texorpdfstring
{Homotopy and face numbers of $X(G_k)$}
{Homotopy and face numbers of the universal complex}}
\label{subsec:universal-complex}

We recall the following facts from \cite{BVV} (see also \cite{Ka}):

\begin{itemize}
    \item[(A)] $X(G_k)$ ($k>1$) is homotopy equivalent to a wedge of $A_{-1,k}$ spheres $S^{k-1}$, where
    \begin{equation}\label{eq:a-k}
    A_{-1,k} = (-1)^k + \sum_{i=0}^{k-1}(-1)^{k-1-i}\frac{\prod_{j=0}^i(2^{k}-2^j)}{(i+1)!}.
    \end{equation}
    \item[(B)] For $-1 \leq p \leq k-2$, $\Lk_{\X}(\sigma^p)$ is homotopy equivalent to a wedge of $A_{p,k}$ spheres $S^{k-p-2}$, where
    \begin{equation}\label{eq:a-pk}
    A_{p,k} = (-1)^{k-p-1} + \sum_{i=0}^{k-p-2} (-1)^{k-p-i} \frac{\prod_{j=1}^{i+1}(2^k - 2^{p+j})}{(i+1)!}.
    \end{equation}
    \item[(C)] The number of $p$-simplices in $\X$ is
    \begin{equation}\label{eq:f-pk}
    f_{p,k} = \frac{(2^k - 2^p) \cdots (2^k - 2^0)}{(p + 1)!}, \quad p \geq 0.
    \end{equation}
\end{itemize}

For completeness, (B) also follows directly from the shelling used
in \cite{BVV}.  Links in a pure shellable complex are shellable;
hence every link here has the homotopy type of a wedge of spheres in
its top dimension $k-p-2$.  The number $A_{p,k}$ is determined by
the Euler characteristic
\[
\chi(\Lk(\sigma^p)) = 1 + (-1)^{k-p-2}A_{p,k} = \sum_{i=0}^{k-p-2} (-1)^i f_i',
\]
where
\[
f_i' = \frac{\prod_{j=1}^{i+1}(2^k - 2^{p+j})}{(i+1)!}
\]
is the number of $i$-simplices in $\Lk(\sigma^p)$.

\subsection{\texorpdfstring
{The spectral sequence of the double complex
$\{E_{p,q},d_h,d_v\}$}
{The spectral sequence of the double complex}}

We refer to \cite[Chapter 10]{Rotman} for the construction of the
spectral sequence associated to a double complex.
We denote the differential on the $E^r$ page by
\[
 d^r\colon E^r_{p,q}\longrightarrow E^r_{p-r,q+r-1}.
\]
Thus $d_h$ and $d_v$ always denote the two double-complex
differentials, whereas $d^r$ is reserved for a spectral-sequence
differential.

We first treat all $k\geq2$ uniformly.  The $E^0$ page has the
following schematic form.

\vskip .1cm
\textbf{$E^0$ Page}:
\[
\begin{array}{|c|cccccc|}
\cdots &  &  & &\cdots &  & \\
k-1&0 &0 &0&\cdots & 0 & 0\\
k-2 & E_{-1,k-2} &E_{0,k-2} & 0&\cdots&0 & 0 \\
k-3 & E_{-1,k-3}&E_{0,k-3} & E_{1,k-3} &\cdots & 0 & 0 \\
\cdots &  &  & &\cdots &  & \\
0 & E_{-1,0} & E_{0,0} & E_{1,0}&\cdots& E_{k-2,0} & 0 \\
-1 & E_{-1,-1}& E_{0,-1} & E_{1,-1} &\cdots& E_{k-2,-1} & E_{k-1,-1} \\ \hline
(p,q) & -1 & 0 & 1 &\cdots & k-2 & k-1 \\ \hline
\end{array}
\]

Taking vertical homology and using the results on the universal complex in
Subsection~\ref{subsec:universal-complex}, the homology of vertical
differentials vanishes in many places:
\begin{equation}\label{eq:e1-page}
E^1_{p,q} =
\begin{cases}
    \bigoplus\limits_{\overline{\sigma^p}} \Span\{c(\sigma^p)\}\otimes\widetilde H_{k-p-2}(\Lk_{\X}(\sigma^p)), & 0\leq p\leq k-1,\ q=k-p-2,\\
   G_k \otimes \ker d^X_{k-2}, & (p,q)=(-1,k-2),\\
   0, & \text{otherwise.} 
\end{cases}
\end{equation}
Here we have used $\ker d^X_{k-2} = \operatorname{im} d^X_{k-1}$, which follows from $\widetilde H_{k-2}(\X)=0$ by (A). In particular, $E^1_{k-1,-1}=\bigoplus\limits_{\overline{\sigma^{k-1}}} \Span\{c(\sigma^{k-1})\}\otimes \Z_2$.

The $E^1$ page has the following schematic form.

\vskip .1cm 
\textbf{$E^1$ Page}:
\[
\begin{array}{|c|cccccc|}
\cdots &  &  & &\cdots &  & \\
k-1&0 &0 &0&\cdots & 0 & 0\\
k-2 & E^1_{-1,k-2} & E^1_{0,k-2} & 0 &\cdots&0 & 0 \\
k-3 & 0 &0 & E^1_{1,k-3} &\cdots & 0 & 0 \\
\cdots &  &  & &\cdots &  & \\
0 & 0 & 0 & 0 &\cdots& E^1_{k-2,0} & 0 \\
-1 & 0 & 0 & 0 &\cdots& 0 &  E^1_{k-1,-1}\\ \hline
(p,q) & -1 & 0 & 1 &\cdots & k-2 & k-1 \\ \hline
\end{array}
\]
For every $k\geq2$, the only possible nonzero $d^1$-differential is
\begin{equation}
    \begin{aligned}
        d^1:& E_{0,k-2}^1 &\longrightarrow &E_{-1,k-2}^1 \\
          & c(\{\alpha\})\otimes z &\longmapsto & \alpha\otimes z,
    \end{aligned}
\end{equation}
where $\alpha\in G_k\setminus\{0\}$ and
$z\in \widetilde H_{k-2}(\Lk(\{\alpha\}))$.

The map $d^1$ need not be surjective.  The next step is to prove that
$d^2$ maps onto $\operatorname{coker}d^1$.  The following chain-level lifting
provides the required preimages.

\begin{lemma}\label{lem:d2-lifting}
    Assume $k\geq2$.  For each generator $\alpha_0\otimes d^X_{k-1}(\sigma) \in E^1_{-1,k-2}$, where $\alpha_0\in G_k$ and $\sigma$ is a $(k-1)$-simplex in $\X$, there exist $x\in E_{0,k-2}$ and $y\in E_{1,k-3}$ such that
    \[
    d_h(x)= \alpha_0\otimes d^X_{k-1}(\sigma),\qquad
    d_v(x)=d_h(y),\qquad d_v(y)=0.
    \]
\end{lemma}

\begin{proof}
  Suppose $\sigma = \{\alpha_1,\dots,\alpha_k\}\in \X$. Then $\sigma$ is a basis of $G_k$, and $\alpha_0$ can be written uniquely as the sum of some subset of $\sigma$.    

    (1) If $\alpha_0 = \alpha_1 + \cdots + \alpha_k$, then all proper faces of $\sigma$ belong to $\Lk(\alpha_0)$. Set 
    \[
    x = c(\{\alpha_0\})\otimes d^X_{k-1}(\sigma)\in E_{0,k-2},\qquad y=0.
    \]
    Then $d_h(x) = \alpha_0\otimes d^X_{k-1}(\sigma)$,
    $d_v(x)=0=d_h(y)$, and $d_v(y)=0$.
    
    (2) Now suppose $\alpha_0 =\alpha_1+\cdots+\alpha_{k-1}$. Set 
    \[
    \eta = \{\alpha_1,\dots,\alpha_{k-1}\}\in \X,\qquad
    z_1 = \eta + d^X_{k-1}(\sigma)\in \widetilde C_{k-2}(\Lk(\alpha_0)).
    \] 
    Then
    \[
    d^X_{k-2}(z_1)=d^X_{k-2}(\eta),
    \]
    and
    \[
    \alpha_0\otimes d^X_{k-1}(\sigma)
    = \alpha_0\otimes \eta + \alpha_0\otimes z_1
    = (\sum_{i=1}^k\alpha_i)\otimes\eta + \alpha_k\otimes \eta + \alpha_0\otimes z_1.
    \]
    Note that $\eta\in \Lk(\sum_{i=1}^k\alpha_i)$ and $\eta\in \Lk(\alpha_k)$.
    Define
    \[
    x = c(\{\sum_{i=1}^k\alpha_i\})\otimes\eta + c(\{\alpha_k\})\otimes \eta + c(\{\alpha_0\})\otimes z_1\in E_{0,k-2},
    \]
    and
    \[
    y = c(\{\alpha_0, \alpha_k\})\otimes d^X_{k-2}(\eta)\in E_{1,k-3}.
    \]
    Then
    \[
    \begin{aligned}
    d_h(x) &= (\sum_{i=1}^k\alpha_i)\otimes\eta + \alpha_k\otimes \eta + \alpha_0\otimes z_1\\
    &= \alpha_0\otimes (\eta+z_1)
    =\alpha_0\otimes d^X_{k-1}(\sigma),
    \end{aligned}
    \]
    and
    \[
    \begin{aligned}
    d_v(x) &= c(\{\sum_{i=1}^k\alpha_i\})\otimes d^X_{k-2}(\eta) + c(\{\alpha_k\})\otimes d^X_{k-2}(\eta) + c(\{\alpha_0\})\otimes d^X_{k-2}(z_1)\\
    &= (c(\{\alpha_0+\alpha_k\}) + c(\{\alpha_k\}) + c(\{\alpha_0\}))\otimes d^X_{k-2}(\eta).
    \end{aligned}
    \]
    On the other hand,
    \[
    d_h(y) = (c(\{\alpha_0+\alpha_k\}) + c(\{\alpha_k\}) + c(\{\alpha_0\}))\otimes d^X_{k-2}(\eta),
    \]
    so $d_v(x)=d_h(y)$. Also $d_v(y)=0$.

    (3) Note that $\alpha_k = (\alpha_1 + \cdots + \alpha_{k-1}) + (\alpha_1 + \cdots + \alpha_k)$ and all differentials are linear. Since the lemma holds for case (1) and case (2), it holds for $\alpha_0=\alpha_k$. By symmetry, it holds for each $\alpha_0=\alpha_i$, $i=1,\dots,k$. Since $\{\alpha_1,\dots,\alpha_k\}$ is a basis of $G_k$, it holds for all $\alpha_0\in G_k$.
\end{proof}

On the $E^2$ page, the terms affected by $d^1$ become
\[
 E^2_{0,k-2}=\ker d^1,
 \qquad
 E^2_{-1,k-2}=\operatorname{coker}d^1,
\]
while the other terms on the diagonal are unchanged.  The only
possible nonzero differential is
$d^2\colon E_{1,k-3}^2\rightarrow E_{-1,k-2}^2$, where
\[
E_{1,k-3}^2 = E_{1,k-3}^1 = \bigoplus_{\overline{\sigma^1}} \Span\{c(\sigma^1)\} \otimes \widetilde H_{k-3}(\Lk(\sigma^1)),
\]
and
\[
E_{-1,k-2}^2 = \operatorname{coker} d^1 = E_{-1,k-2}^1/\operatorname{im}(d^1).
\]

Explicitly, if $\bar y\in E^2_{1,k-3}$ is represented by a vertical
cycle $y$ and $x\in E_{0,k-2}$ satisfies $d_v(x)=d_h(y)$, then
\[
 d^2(\bar y)=[d_h(x)]\in E^2_{-1,k-2}=\operatorname{coker}d^1.
\]

\begin{proposition}\label{prop:d2-surjective}
    For $k\geq2$, the map $d^2: E_{1,k-3}^2\rightarrow E_{-1,k-2}^2$ is surjective. Consequently, the spectral sequence degenerates at the $E^3$ page, i.e. $E^3 = E^\infty$.
\end{proposition}

\begin{proof}
Every class in $E^2_{-1,k-2}=\operatorname{coker}d^1$ is generated
by the image of some
$\alpha_0\otimes d^X_{k-1}(\sigma)$.  For the elements $x$ and $y$
provided by Lemma~\ref{lem:d2-lifting}, the defining formula for
$d^2$ gives
\[
 d^2([y])=[d_h(x)]
 =[\alpha_0\otimes d^X_{k-1}(\sigma)].
\]
Thus $d^2$ is surjective.  All terms on the $E^2$ page lie on the
total-degree-$(k-2)$ diagonal except its target at $(-1,k-2)$, which
is killed by $d^2$; degree considerations then force all later
differentials to vanish.
\end{proof}

Consequently, the $E^3=E^\infty$ page is supported entirely on the
diagonal $p+q=k-2$.  Its possible nonzero terms are $\ker d^1$ at
$(0,k-2)$, $\ker d^2$ at $(1,k-3)$, and
$E^1_{p,k-2-p}$ for $2\le p\le k-1$.

It remains to consider $k=1$.  If $\alpha_1$ is the nonzero element
of $G_1$, the $E^1$ page has only the two terms
\[
 E^1_{0,-1}\cong\mathbb Z_2,
 \qquad E^1_{-1,-1}\cong\mathbb Z_2,
\]
and
\[
 d^1\colon E^1_{0,-1}\longrightarrow E^1_{-1,-1},
 \qquad
 c(\{\alpha_1\})\otimes\emptyset
 \longmapsto \alpha_1\otimes\emptyset,
\]
is an isomorphism.  Thus $E^2=E^\infty=0$, and
$\B_\bullet(G_1)$ is also acyclic.

\begin{corollary}[Homology concentration]
\label{cor:homology-concentration}
The complex $\B_\bullet(G_k)$ has trivial homology outside
degree $k-2$:
\[
 H_i\bigl(\B_\bullet(G_k)\bigr)=0
 \qquad(i\ne k-2).
\]
\end{corollary}

\begin{proof}
For $k=1$, the assertion follows from the separate computation above.
Now let $k\geq2$.
The bounded spectral sequence converges to the homology of
$\B_\bullet(G_k)$.  Since its $E^\infty$ page is supported
in total degree $k-2$, all other homology groups vanish.
\end{proof}

\begin{proposition}\label{prop:total-homology}
    For $k\geq2$, in degree $k-2$ there is a vector-space isomorphism
    \[
    H_{k-2}\bigl(\B_\bullet(G_k)\bigr)
    \cong \bigoplus_{p=0}^{k-1} E_{p,k-2-p}^3
    \quad\text{as a $\mathbb Z_2$-vector space},
    \]
    where for $p=0,\dots,k-1$,
    \[
    E_{p,k-2-p}^3 =
    \begin{cases}
        \ker d^1, & p=0,\\
        \ker d^2, & p=1,\\
        E_{p,k-2-p}^1, & 2\leq p\leq k-1.
    \end{cases}
    \]
\end{proposition}

\begin{proof}
In total degree $k-2$, the $E^\infty$ page gives the associated
graded vector space of a finite filtration of the homology.  Since
the ground ring is the field $\mathbb Z_2$, this filtration splits
(not canonically), giving the displayed vector-space isomorphism.
\end{proof}

\begin{proof}[Proof of Theorem~\ref{thm:intro-dimension}]
For $k=1$, the spectral sequence directly gives
$\mathcal Z_2(G_1)=0$; this case is not covered by the theorem.
Now assume that $k\ge2$.  By
Proposition~\ref{prop:total-homology}, the desired dimension is the
sum of the dimensions of the surviving terms on the
total-degree-$(k-2)$ diagonal.  Set $A_{k-1,k}:=1$, corresponding to
the reduced homology of the link of a top-dimensional simplex.  Then
\eqref{eq:e1-page} and Proposition~\ref{prop:equiv-facts}(4) give
\[
 \dim E^1_{p,k-2-p}
 =\frac{A_{p,k}f_{p,k}}{p+2}
 \qquad(1\le p\le k-1).
\]
Rank--nullity and the surjectivity
\[
 d^2\colon E^2_{1,k-3}\longrightarrow
 E^2_{-1,k-2}=\operatorname{coker}d^1
\]
give the calculation without determining either kernel separately:
\begin{align*}
 \dim\ker d^1+\dim\ker d^2
 &=(\dim E^1_{0,k-2}-\rank d^1)
   +(\dim E^1_{1,k-3}-\dim\operatorname{coker}d^1)\\
 &=\dim E^1_{0,k-2}+\dim E^1_{1,k-3}
   -\dim E^1_{-1,k-2}.
\end{align*}
By \eqref{eq:e1-page}, this equals
\[
 A_{0,k}f_{0,k}
 +\frac{A_{1,k}f_{1,k}}{3}
 -k\bigl(f_{k-1,k}-A_{-1,k}\bigr).
\]
Adding the remaining terms gives
\[
\begin{aligned}
 \dim_{\mathbb Z_2}\mathcal Z_{k+1}(G_k)
 &=A_{0,k}f_{0,k}
   +\sum_{p=1}^{k-1}\frac{A_{p,k}f_{p,k}}{p+2}
   -k\bigl(f_{k-1,k}-A_{-1,k}\bigr)\\
 &=A_{0,k}f_{0,k}
   +\sum_{p=1}^{k-2}\frac{A_{p,k}f_{p,k}}{p+2}
   +kA_{-1,k}
   -\left(k-\frac1{k+1}\right)f_{k-1,k},
\end{aligned}
\]
which is the stated formula.

\end{proof}

\section{Appendix: the horizontal differential}\label{sec:appendix}

We record the two facts about the horizontal differential used in
Proposition~\ref{prop:double-complex}.

\begin{lemma}\label{lem:dh-equals-dC}
For $p\geq2$,
\[
 (d_h)_{p,q}(c(\sigma^p)\otimes\sigma^q)
 =d^X_p(c(\sigma^p))\otimes\sigma^q.
\]
\end{lemma}

\begin{proof}
Write $\sigma^p=\{\alpha_1,\dots,\alpha_{p+1}\}$,
$\sigma_i=\sigma^p\setminus\{\alpha_i\}$, and
$\sigma_{ij}=\sigma^p\setminus\{\alpha_i,\alpha_j\}$.  Applying
$d^X_p$ to \eqref{eq:equiv-sum} gives
\[
\begin{aligned}
 d^X_p(c(\sigma^p))
 ={}&\sum_i\sigma_i+\sum_i(\sigma_i+\alpha_i)\\
 &+\sum_i\sum_{j\ne i}
   ((\sigma_{ij}+\alpha_i)\cup\{\alpha_i\}).
\end{aligned}
\]
On the other hand, \eqref{eq:equiv-sum} gives
\[
 \sum_i c(\sigma_i)
 =\sum_i\sigma_i
  +\sum_i\sum_{j\ne i}
    ((\sigma_{ij}+\alpha_j)\cup\{\alpha_j\}).
\]
Interchanging $i$ and $j$ in the double sum replaces $\alpha_j$ by
$\alpha_i$.  Moreover, \eqref{eq:equiv-sum} and
Lemma~\ref{lem:shift-equivalence} give
\[
 c(\sigma_1+\alpha_1)=\sum_i(\sigma_i+\alpha_i).
\]
Comparing this expression with the defining formula for $(d_h)_{p,q}$
applied to $c(\sigma^p)\otimes\sigma^q$ proves the claimed identity.
\end{proof}

\begin{lemma}\label{lem:dh-well-defined}
For every $p$-simplex $\sigma^p\in X(G_k)$,
$(d_h)_{p,q}(c(\sigma^p)\otimes\sigma^q)$ depends only on
$\overline{\sigma^p}$.
\end{lemma}

\begin{proof}
For $p=0$, the equivalence class is a singleton.  For $p=1$, write
$\sigma^1=\{\alpha_1,\alpha_2\}$.  Set 
$
 W=\Span_{\mathbb Z_2}\sigma^1
  =\{0,\alpha_1,\alpha_2,\alpha_1+\alpha_2\}.
$ Then
\[
 (d_h)_{1,q}(c(\sigma^1)\otimes\sigma^q)
 =\left(\sum_{\alpha\in W\setminus\{0\}}c(\{\alpha\})\right)
  \otimes\sigma^q,
\]
which depends only on $W$, and hence only on $\overline{\sigma^1}$.
For $p\geq2$, the assertion follows from
Lemma~\ref{lem:dh-equals-dC}, since the equivalence chain
$c(\sigma^p)$ depends only on $\overline{\sigma^p}$.
\end{proof}

\begin{proposition}\label{prop:dh-square-zero}
$(d_h)_{p-1,q}\circ(d_h)_{p,q}=0$ for all $p,q$.
\end{proposition}

\begin{proof}
For $p\leq0$, the composite has zero target.  For $p\geq3$,
Lemma~\ref{lem:dh-equals-dC} gives
$d_h^2=(d^X)^2\otimes\mathrm{id}=0$.

For $p=1$, let $\sigma^1=\{\alpha_1,\alpha_2\}$.  Then
\[
 (d_h)_{0,q}\circ(d_h)_{1,q}(c(\sigma^1)\otimes\sigma^q)
 =(\alpha_1+\alpha_2+(\alpha_1+\alpha_2))\otimes\sigma^q=0.
\]

For $p=2$, let $\sigma^2=\{\alpha_1,\alpha_2,\alpha_3\}$ and set
$\alpha_{ij}=\alpha_i+\alpha_j$.  A direct
calculation gives
\[
\begin{aligned}
 (d_h)_{1,q}\circ(d_h)_{2,q}(c(\sigma^2)\otimes\sigma^q)
 ={}&\bigl(
 c(\{\alpha_{12}\})+c(\{\alpha_{13}\})+c(\{\alpha_{23}\})\\
 &\quad+c(\{\alpha_1\})+c(\{\alpha_2\})+c(\{\alpha_{12}\})\\
 &\quad+c(\{\alpha_1\})+c(\{\alpha_3\})+c(\{\alpha_{13}\})\\
 &\quad+c(\{\alpha_2\})+c(\{\alpha_3\})+c(\{\alpha_{23}\})
 \bigr)\otimes\sigma^q=0.
\end{aligned}
\]
\end{proof}

\bibliographystyle{abbrv}

\end{document}